\documentclass{article}

\usepackage{amssymb}
\usepackage{amsmath}

\usepackage[dvips]{graphicx}

\begin{document}


\begin{center}
\textbf{On Saturated Uniformly A-convex Algebras} 
\end{center}

\noindent\textbf{•}

\begin{center}
M. El Azhari
\end{center}
 
\noindent\textbf{•}

\noindent\textbf{•}Department of Mathematics, Ecole Normale Sup\'{e}rieure, Rabat, Morocco.

\noindent\textbf{•}E-mail: mohammed.elazhari@yahoo.fr 

\noindent\textbf{•}

\noindent\textbf{•}

\noindent\textbf{Abstract. }Following ideas of A. C. Cochran, we give a suitable definition of a saturated uniformly A-convex algebra.  In the m-convex case, such algebra is a uniform topological one.

\noindent\textbf{•}

\noindent\textbf{Keywords.} Uniformly A-convex Algebra; Saturated Algebra; Uniform Topological Algebra.

\noindent\textbf{•}

\noindent\textbf{•}

\noindent\textbf{1. PRELIMINARIES}

\noindent\textbf{•}

A locally convex algebra is a complex algebra with a locally convex topology for which the multiplication is separately continuous.  A locally m-convex algebra is a locally convex algebra whose topology is defined by a family of submultiplicative seminorms.  A uniform seminorm on an algebra  $E$ is a seminorm $p$ such that  $p\left(x^2\right)=p(x)^2 $ for all $x\in E.$  Such a seminorm is submultiplicative.  A uniform topological algebra is a locally convex algebra whose topology is defined by a family of uniform seminorms.  A uniform normed algebra is a normed algebra $(E,\Vert .\Vert)$  such that $ \Vert x^{2}\Vert=\Vert x\Vert^{2}$ for all $x\in E. $ A locally convex algebra $E$ is uniformly A-convex if its topology is defined by a family $\lbrace p_{\alpha},\alpha\in \Lambda\rbrace$ of seminorms with the property that for $x\in E,$  there is a positive constant  $r_x$ such that $p_{\alpha }(xy)\le r_x p_{\alpha }(y)$  and  $p_{\alpha }(yx)\le r_x p_{\alpha }(y)$  for all $\alpha \in \Lambda $  and  $y\in E.$  Let $E$ be a locally convex algebra. Denote by  $M^{\ast}\left(E\right)$  the set of all nonzero multiplicative linear functionals on $E$. Denote by   $M(E)$  the space of all nonzero continuous multiplicative linear functionals on $E,$ topologized via the weak topology,  it is called the carrier space of $E.$

\noindent\textbf{•}

\noindent\textbf{2. RESULTS}

\noindent\textbf{•}

Let $(E,\left(p_{\alpha }\right)_{\alpha \in \Lambda})$ be a  Hausdorff commutative uniformly A-convex algebra  with unit $e.$  The family $\lbrace p_{\alpha},\alpha \in \Lambda\rbrace$ of seminorms can be chosen such that  $p_{\alpha }\left(e\right)=1 $ for all $\alpha \in \Lambda.$  For $x\in E,$  let $\Vert x\Vert=\sup  \left[\sup\lbrace p_{\alpha }\left(xy\right), p_{\alpha }\left(y\right)\le 1\rbrace:\alpha \in \Lambda \right]=\inf \lbrace r_x>0, p_{\alpha }\left(xy\right)\le r_{x}p_{\alpha }\left(y\right) $ for all $ \alpha \in \Lambda, y\in E\rbrace . $ By [1, Lemma 3.2],  $ \Vert .\Vert $  is a submultiplicative norm on $E$ for which $ p_{\alpha }(x)\le \Vert x\Vert $  for all $\alpha \in \Lambda $ and $x\in E.$  Let $M$ be the carrier space of $(E,\left(p_{\alpha }\right)_{\alpha \in \Lambda })$ and let $M_n$ be the carrier space of $\left(E,\Vert .\Vert\right), M\subset M_n$.  Since $\left(E,\Vert .\Vert\right)$ is a commutative normed algebra with unit, $M_n$ is nonempty. $M$ may be empty [2].  For the sequel, we assume that $M$ is nonempty.  For $\alpha \in \Lambda $ and $m\in M,$ Cochran [1] has defined the extended real number  $t_{\alpha }\left(m\right)=\sup\lbrace\vert m(x)\vert ,p_{\alpha }(x)\le 1\rbrace $  and the map $\phi_{\alpha }:M\to R $ by $\phi_{\alpha }\left(m\right)=t_{\alpha }(m)^{-1}$  if $t_{\alpha }(m)<\infty $  and $\phi_{\alpha }\left(m\right)=0$  otherwise.

\noindent\textbf{•}

\noindent\textbf{Proposition 2.1.} 
 
(1)\quad $t_{\alpha }(m)\ge 1$ and  $0\le \phi_{\alpha }(m)\le 1$  for all $\alpha \in \Lambda $  and $m\in M;$   

(2)\quad $t_{\alpha }(m)<\infty $  if, and only if, $m$ is continuous for $p_{\alpha };$           

(3)\quad If $t_{\alpha }\left(m\right)<\infty ,$ then $\vert m(x)\vert\le t_{\alpha }\left(m\right) p_{\alpha }(x)$  for all $x\in E.$   
 
\noindent\textbf{•}

\noindent\textbf{Proof.} (1)\quad  Since $m\left(e\right)=1$ and $p_{\alpha }\left(e\right)=1,$ it follows that $t_{\alpha }(m)\ge 1,$ hence $\phi_{\alpha }\left(m\right)=t_{\alpha }(m)^{-1}\le 1.$

(2)\quad $m$ is continuous for $p_{\alpha }$ if, and only if, $m$ is bounded on $\lbrace x\in E, p_{\alpha }(x)\le 1\rbrace$, i.e.  $t_{\alpha }\left(m\right)<\infty .$

(3)\quad Let $x\in E$ and $\varepsilon >0,\vert m\left(\left(p_{\alpha }\left(x\right)+\varepsilon \right)^{-1}x\right)\vert\le t_{\alpha }(m)$  since $ p_{\alpha }\left(\left(p_{\alpha }\left(x\right)+\varepsilon \right)^{-1} x)\right)\le 1,$ thus  $\vert m\left(x\right)\vert\le \ t_{\alpha }\left(m\right)\left(p_{\alpha }\left(x\right)+\varepsilon \right).$ Since $\varepsilon >0$ is arbitrary, we conclude that $\vert m\left(x\right)\vert\leq t_{\alpha }\left(m\right)p_{\alpha }\left(x\right).$

\noindent\textbf{•}
 
Let $C_{b}(M)$ be the algebra of all complex continuous bounded functions on $M,$ with the topology defined by the family $\lbrace\hat{p}_{\alpha } ,\alpha \in \Lambda\rbrace$ of seminorms, where $\hat{p}_{\alpha}(f)=\sup \lbrace \phi_{\alpha}(m)\vert f(m)\vert, m\in M\rbrace $ for all $f\in C_{b}\left(M\right). $ For $\alpha \in \Lambda $ , define $M_{\alpha }=\lbrace m\in M, t_{\alpha }\left(m\right)<\infty\rbrace.$ If $M_{\alpha }$ is empty, $\phi_{\alpha }\left(m\right)=0$ for all $m\in M,$ so $\hat{p}_{\alpha }\left(f\right)=0$ for all $f\in C_{b}\left(M\right).$ If $M_{\alpha }$ is nonempty,  $\hat{p}_{\alpha }(f)=\sup \lbrace t_{\alpha}(m)^{-1}\vert f(m)\vert, m\in M_{\alpha}\rbrace .$ If $\lbrace p_{\alpha } , \alpha \in \Lambda \rbrace$ is a directed family of seminorms, then $M=\bigcup_{\alpha \in \Lambda }M_{\alpha }$ , so there exists $\alpha \in \Lambda $ such that $M_{\alpha }$ is nonempty. The following result is due to Cochran [1], the proof is given for completeness.

\noindent\textbf{•}

\noindent\textbf{Proposition 2.2 [1].} Let $G:E\to C(M)$ be the Gelfand map. Then
 
(1)\quad $ \vert\hat{x}(m)\vert\le \Vert x\Vert$  for all $x\in E$ and $m\in M,$ so $G\left(E\right)=\hat{E}\subset C_{b}\left(M\right);$  

(2)\quad $ \hat{p}_{\alpha}\left(\hat{x}\right)\le p_{\alpha }(x)$ for all $\alpha \in \Lambda $  and $x\in E.$

\noindent\textbf{•}

\noindent\textbf{Proof.} (1)\quad Let  $m\in M\subset M_n$. Since $(E,\Vert .\Vert)$ is a normed algebra and  $m\in M_n$ , it follows that  $\vert\hat{x}\left(m\right)\vert=\vert m(x)\vert\le \Vert x\Vert $  for all $x\in E.$ Consequently ,  $\hat{E}\subset C_b\left(M\right).$

(2)\quad Let $\alpha \in \Lambda.$ If  $M_{\alpha }$ is empty,  ${\hat{p}}_{\alpha }\left(\hat{x}\right)=0\le p_{\alpha }(x)$  for all $x\in E.$ If $M_{\alpha }$ is nonempty, we have $\vert m\left(x\right)\vert\le t_{\alpha }\left(m\right) p_{\alpha }(x)$  for all $ m\in M_{\alpha }$ and $x\in E,$ i.e. $ t_{\alpha }\left(m\right)^{-1}\vert m\left(x\right)\vert\le 
p_{\alpha }(x)$ for all $m\in M_{\alpha }$ and $x\in E.$ Thus  $\hat{p}_{\alpha }\left(\hat{x}\right)= \sup \lbrace t_{\alpha }\left(m\right)^{-1}\vert m\left(x\right)\vert, m\in M_{\alpha }\rbrace\le p_{\alpha }(x)$ for all $x\in E.$

\noindent\textbf{•}

Let $\alpha \in \Lambda $  and let $F_{\alpha }$ be the topological dual of $\left(E,p_{\alpha }\right). $ For $f\in F_{\alpha }$, put $ t_{\alpha }\left(f\right)=\sup\lbrace\vert f(x)\vert, p_{\alpha}(x)\leq 1\rbrace, t_{\alpha }$ is a norm on $F_{\alpha }$. If $M_{\alpha }$ is nonempty, $M_{\alpha }$ is usually topologized via the weak topology. But for the following proposition, we consider on $M_{\alpha }$ the topology $\tau_{\alpha }$ induced by the topology of $t_{\alpha }$. The map $M_{\alpha }\to R, m\to t_{\alpha }\left(m\right),$ is continuous for $\tau_{\alpha }$. Since $(M_{\alpha } ,\tau_{\alpha })$ is a metric space, it is a completely regular space, we denote by $M^{\beta }_{\alpha }$ the Stone-$\check{C}$ech compactification of $\ M_{\alpha }$ .

\noindent\textbf{•}

\noindent\textbf{Proposition 2.3.} The following assertions are equivalent:

(1)\quad  If $ p_{\alpha }\left(x\right)=1 $ for some $\alpha \in \Lambda $ and $x\in E,$ then  $\hat{p}_{\alpha }\left(\hat{x}\right)=1\ ;$

(2)\quad  $\hat{p}_{\alpha }\left(\hat{x}\right)=p_{\alpha }\left(x\right)\ $ for all $\alpha \in \Lambda $  and $x\in E ;$

(3)\quad  For  $\alpha \in \Lambda $  and $x\in E $ such that  $p_{\alpha }\left(x\right)=1 , M_{\alpha }$ is nonempty and the continuous extension of the map  $M_{\alpha }\to R, m\to t_{\alpha }\left(m\right)^{-1}\vert m\left(x\right)\vert,$  to $M^{\beta }_{\alpha }$ is equal to 1 at some $ m_{0}\in M^{\beta }_{\alpha }$.

\noindent\textbf{•}

\noindent\textbf{Proof.} (1) $\Rightarrow $ (2)\quad : Let $ \alpha \in \Lambda  $ and $ x\in E$. If $p_{\alpha }\left(x\right)=0,$ then  $0\le p_{\alpha }\left(\hat{x}\right)\le p_{\alpha }\left(x\right)=0 $ by Proposition 2.2, so  $\hat{p}_{\alpha }\left(\hat{x}\right)=0.$ If $ p_{\alpha }\left(x\right)\ne 0,\hat{p}_{\alpha }\left(p_{\alpha }\left(x\right)^{-1} \hat{x}\right)=1$  since $p_{\alpha }\left(p_{\alpha }\left(x\right)^{-1} x\right)=1, $ hence  $\hat{p}_{\alpha }\left(\hat{x}\right)=p_{\alpha }\left(x\right).$

(2) $\Rightarrow $ (3)\quad :  Let  $\alpha \in \Lambda $  and $x\in E $ such that  $p_{\alpha }\left(x\right)=1.$ $M_{\alpha }$ is nonempty since  $1=p_{\alpha }\left(x\right)=\hat{p}_{\alpha }\left(\hat{x}\right)\ne 0 $.  Let $m_1, m_2\ $ in $M_{\alpha }$ , $\vert m_1\left(x\right)-m_2\left(x\right)\vert=\vert\left(m_1-m_2\right)\left(x\right)\vert\le t_{\alpha }\left(m_1-m_2\right)p_{\alpha }\left(x\right)=t_{\alpha }\left(m_1-m_2\right),$  then the map $M_{\alpha }\to R, m\to m\left(x\right), $ is continuous. Consequently, the map $\varphi_{\alpha }: M_{\alpha }\to R, \varphi _{\alpha }\left(m\right)=t_{\alpha }\left(m\right)^{-1}\vert m\left(x\right)\vert,$  is continuous.  $\varphi_{\alpha }$ is also bounded on $M_{\alpha }$ since $\hat{p}_{\alpha }\left(\hat{x}\right)=\sup\lbrace t_{\alpha}(m)^{-1}\vert m(x)\vert, m\in M_{\alpha }\rbrace=p_{\alpha }\left(x\right)=1.$  Therefore  $\varphi_{\alpha }$ has a continuous extension $\varphi^{\beta }_{\alpha }$ to $M^{\beta }_{\alpha }$. Since $ t_{\alpha }\left(m\right)^{-1}\vert m\left(x\right)\vert\le 1$ for all $ m\in M_{\alpha }$ and $M_{\alpha }$ is dense in  $M^{\beta }_{\alpha }$ , it follows that $\varphi ^{\beta }_{\alpha }\left(m\right)\le 1$ for all $m\in M^{\beta }_{\alpha }$.  We have $ 1= \sup \lbrace t_{\alpha }\left(m\right)^{-1}\vert m(x)\vert, m\in M_{\alpha }\rbrace\le \sup \lbrace\varphi ^{\beta }_{\alpha }\left(m\right), m\in M^{\beta }_{\alpha }\rbrace\le 1,$  hence $\sup \lbrace  \varphi ^{\beta }_{\alpha }\left(m\right), m\in M^{\beta }_{\alpha }\rbrace =1$. Since $M^{\beta }_{\alpha }$ is compact, there exists $m_0\in M^{\beta }_{\alpha }$  such that $\varphi^{\beta }_{\alpha }\left( m_0\right)=1.$

(3)$ \Rightarrow $ (1)\quad :  Let  $\alpha \in \Lambda $  and $x\in E $ such that   $p_{\alpha }\left(x\right)=1. $ If $\hat{p}_{\alpha }\left(\hat{x}\right)=\sup\lbrace t_{\alpha }\left(m\right)^{-1}\vert m(x)\vert, m\in M_{\alpha}\rbrace=s<1, t_{\alpha}(m)^{-1}\vert m(x)\vert\leq s<1$ for all $ m\in M_{\alpha} $ , then  $\varphi^{\beta }_{\alpha }\left(m\right)\le s<1$ for all $m\in M^{\beta }_{\alpha }$, this contradicts (3)\quad .  Therefore $\hat{p}_{\alpha }\left(\hat{x}\right)\ge 1,$ so $1\le \hat{p}_{\alpha }\left(\hat{x}\right)\le p_{\alpha }\left(x\right)=1 $  i.e. $ \hat{p}_{\alpha }\left(\hat{x}\right)=1.$

\noindent\textbf{•}

We say that $E$ is a Cochran algebra if $\hat{p}_{\alpha }\left(\hat{x}\right)=p_{\alpha }\left(x\right)$ for all $\alpha \in \Lambda $  and $x\in E.$
 
\noindent\textbf{•}

\noindent\textbf{Proposition 2.4.} If $E$ is a Cochran algebra, then the Gelfand map $G$ from $(E,\left(p_{\alpha }\right)_{\alpha \in \Lambda })$  onto $(\hat{E},\left(\hat{p}_{\alpha }\right)_{\alpha \in \Lambda })$ is an algebraic and topological isomorphism.

\noindent\textbf{•}

\noindent\textbf{Proof.} Let $ x\in E.$ If $\hat{x}=0, $ then  $\hat{p}_{\alpha }\left(\hat{x}\right)=p_{\alpha }\left(x\right)=0 $ for all $\alpha \in \Lambda  , $ so $x=0$ since $E$ is Hausdorff. Consequently, $G$ is an algebraic isomorphism. Since  $\hat{p}_{\alpha }\left(\hat{x}\right)=p_{\alpha }\left(x\right) $ for all $\alpha \in \Lambda $  and $x\in E,$  $G$ is a topological isomorphism.

\noindent\textbf{•}

\noindent\textbf{Proposition 2.5.} Let $E$ be a Cochran algebra whose carrier space $M$ is equicontinuous. Then $E$ is a uniform normed algebra.

\noindent\textbf{•}

\noindent\textbf{Proof. }For $x\in E,$ let $q\left(x\right)=\sup \lbrace \vert m(x)\vert, m\in M\rbrace.$  Since $M$ is equicontinuous, the map $q$ is a continuous seminorm on $E.$ Let  $\alpha \in \Lambda $  and $x\in E, p_{\alpha }\left(x\right)=\hat{p}_{\alpha }\left(\hat{x}\right)=\sup \lbrace \phi_{\alpha}(m)\vert m(x)\vert, m\in M\rbrace\leq q(x)$ since $0\leq \phi_{\alpha}(m)\leq 1  $  for all $m\in M.$ Let $x\in E$ such that $q\left(x\right)=0,$  since $E$ is Hausdorff and $p_{\alpha }\left(x\right)\le q\left(x\right)=0$  for all $\alpha \in \Lambda $ , it follows that $x=0.$  $q$ is a continuous uniform norm on $E$ such that $ p_{\alpha }(x)\le q\left(x\right)$  for all $\alpha \in \Lambda $  and $x\in E,$ then the topology of $E$ can be defined by the uniform norm $q.$

\noindent\textbf{•}

\noindent\textbf{Corollary 2.1.} Let $E$ be a Cochran algebra. If $E$ is barrelled or a Q-algebra, then $E$ is a uniform normed algebra.

\noindent\textbf{•}
 
\noindent\textbf{Proof.} (1)\quad  Assume that $E$ is barrelled. By Proposition 2.2, $\vert m(x)\vert\le \Vert x\Vert$  for all $x\in E$ and $m\in M.$ Therefore $M$ is bounded for the weak topology, so $M$ is equicontinuous.

(2)\quad  Assume that $E$ is a Q-algebra. By [3, Proposition II.7.1], every topological Q-algebra has an equicontinuous carrier space.

\noindent\textbf{•}

The saturated uniformly A-convex algebras were introduced by Cochran [1] as a subclass of the class of Cochran algebras.

\noindent\textbf{•}

\noindent\textbf{Definition 2.1 [1].} $E$ is saturated if for each $ \alpha \in \Lambda $ and $x\in E$ such that $p_{\alpha} \left(x\right)=1$, there exists $ m_0\in M$  such that $ m_0\left(x\right)=\sup \lbrace\vert m(y)\vert, p_{\alpha }(y)\le 1\rbrace$ for some $m\in M.$ ( Then $\hat{p}_{\alpha }\left(\hat{x}\right)=1=p_{\alpha }\left(x\right)$).

\noindent\textbf{•}

Oudadess [2] has proved that the class of complete saturated algebras  (in the sense of Definition 2.1) is empty. His proof uses the completeness of the algebra.  Here we give a simple proof in the general case.

\noindent\textbf{•}

\noindent\textbf{Proposition 2.6.} The class of saturated algebras (in the sense of Definition 2.1) is empty.

\noindent\textbf{•}
 
\noindent\textbf{Proof.} Let $ \alpha \in \Lambda $ , we have  $p_{\alpha }\left(-e\right)=p_{\alpha }\left(e\right)=1.$  Since the algebra is saturated, there exist $ m_0, m$ in $M$ such that $\ -1=m_0\left(-e\right)=\sup\lbrace\vert m(y)\vert, p_{\alpha }(y)\le 1\rbrace\geq 0,$  which is absurd.

\noindent\textbf{•}

We think that the definition of a saturated algebra should be as follows:

\noindent\textbf{•}

\noindent\textbf{Definition 2.2.} $E$ is saturated if for each $\alpha \in \Lambda $  and $x\in E$ such that  $p_{\alpha }\left(x\right)=1,$ there exists  $m_0\in M$  such that $\vert m_0\left(x\right)\vert=\sup \lbrace\vert m_{0}(y)\vert, p_{\alpha}(y)\leq 1\rbrace,$ i.e. $ m_{0}\in M_{\alpha}$ and $ t_{\alpha}(m_{0})^{-1}\vert m_{0}(x)\vert=1. $

\noindent\textbf{•}

\noindent\textbf{Proposition 2.7.} If $E$ is saturated, then $E$ is a Cochran algebra.

\noindent\textbf{•}
 
\noindent\textbf{Proof.} Let $\alpha \in \Lambda $  and $x\in E $ such that  $p_{\alpha }\left(x\right)=1.$ By hypothesis, $M_{\alpha }$ is nonempty and the map  $M_{\alpha }\to R, m\to t_{\alpha }\left(m\right)^{-1}\vert m\left(x\right)\vert,$ is equal to 1 at some $m_0\in M.$  Therefore  $E$ is a Cochran algebra by Proposition 2.3.

\noindent\textbf{•}
 
If $p_{\alpha }$ is submultiplicative, $N_{\alpha }=\lbrace x\in E, p_{\alpha }\left(x\right)=0\rbrace$ is an ideal in $E$ and the quotient algebra $ E/N_{\alpha }$ is a normed algebra with the norm $ \Vert x_{\alpha}\Vert_{\alpha}=p_{\alpha }\left(x\right), x_{\alpha }=x+N_{\alpha } .$ Let $E_{\alpha }$ be the completion of $ E/N_{\alpha },$ $ E_{\alpha }$ is a commutative Banach algebra with unit.

\noindent\textbf{•}

\noindent\textbf{Proposition 2.8.} If $p_{\alpha }$ is submultiplicative, then

(1)\quad  $M_{\alpha }$ is a nonempty compact space;

(2)\quad  $ t_{\alpha }\left(m\right)=1$  for all $m\in M_{\alpha }.$

\noindent\textbf{•}
 
\noindent\textbf{Proof.} (1)\quad  $M_{\alpha }$ is homeomorphic to $M(E_{\alpha })$ by [4, Proposition 7.5], then $M_{\alpha }$ is a nonempty compact space. 

(2)\quad  Let $m\in M_{\alpha },\vert m\left(x\right)\vert\le t_{\alpha }\left(m\right) p_{\alpha }(x)$  for all $x\in E,$ then $N_{\alpha }$ is included in $Ker\left(m\right).$ We may define a multiplicative linear functional $m_{\alpha }$ on $E/N_{\alpha }$ by $m_{\alpha }\left(x_{\alpha }\right)=m\left(x\right).$  Since $E/N_{\alpha }$ is a normed algebra and $m_{\alpha }\in M(E/N_{\alpha }),\vert m_{\alpha }\left(x_{\alpha }\right)\vert\le \Vert x_{\alpha}\Vert_{\alpha}$ for all $x\in E.$  Further, as $ E/N_{\alpha }$ is unital, $ 1=\sup \lbrace\vert m_{\alpha}(x_{\alpha})\vert, \Vert x_{\alpha}\Vert_{\alpha}\leq 1\rbrace=\sup \lbrace\vert m(x)\vert, p_{\alpha}(x)\leq 1\rbrace=t_{\alpha}(m).$

\noindent\textbf{•}
 
\noindent\textbf{Proposition 2.9.} If $ p_{\alpha }$ is submultiplicative for every $\alpha \in \Lambda $, then the following assertions are equivalent:

(1)\quad $E$ is saturated;

(2)\quad $E$ is a Cochran algebra;

(3)\quad $(E,\left(p_{\alpha }\right)_{\alpha \in \Lambda })$  is a uniform topological algebra.

\noindent\textbf{•}

\noindent\textbf{Proof.} (1) $\Rightarrow $(2)\quad: By Proposition 2.7.

(2) $\Rightarrow $ (3)\quad:   Let  $\alpha \in \Lambda $  and $x\in E$,  $p_{\alpha }\left(x\right)=\hat{p}_{\alpha }\left(\hat{x}\right)=\sup \lbrace t_{\alpha}(m)^{-1}\vert m(x)\vert, m\in M_{\alpha}\rbrace=\sup \lbrace \vert m(x)\vert, m\in M_{\alpha}\rbrace$ by Proposition 2.8. Thus  $p_{\alpha }$ is  a uniform seminorm.

(3) $\Rightarrow $ (1)\quad:  Let  $\alpha \in \Lambda $  and $x\in E $ such that  $p_{\alpha }\left(x\right)=1.$  Since $E_{\alpha }$ is a uniform Banach algebra, $ 1= p_{\alpha}(x)=\Vert x_{\alpha}\Vert_{\alpha}=\sup \lbrace\vert g(x_{\alpha})\vert, g\in M(E_{\alpha})\rbrace=\sup \lbrace\vert m(x)\vert, m\in M_{\alpha}\rbrace$ by [4, proposition 7.5].  Since $M_{\alpha }$ is compact, there exists  $m_0\in M_{\alpha }$ such that $\vert m_0\left(x\right)\vert=1=t_{\alpha }\left(m_0\right)$  by Proposition 2.8, hence $E$ is saturated.

\noindent\textbf{•}

We give an example of a saturated algebra which is not m-convex.

\noindent\textbf{•}
 
\noindent\textbf{Example.} Let $C_b(R)$ be the algebra of all complex continuous bounded functions on $R.$ Let $ \Lambda $ be the subset of $C_b(R)$ defined by  $\phi \in \Lambda $ if $\phi(x)>0$ for all $x\in R, \sup\lbrace\vert \phi (x)\vert, x\in R\rbrace =1 $  and $ \phi$  vanishes at infinity. We endow $C_b(R)$ with the topology determined by the family $\lbrace p_{\phi} , \phi\in \Lambda\rbrace $ of seminorms, where $ p_{\phi}(f)=\sup \lbrace\vert f(x)\phi (x)\vert, x\in R\rbrace.$ $C_b(R)$ is a commutative complete uniformly A-convex algebra with unit $u (u\left(x\right)=1$ for all $x\in R)$ such that $ p_{\phi}\left(u\right)=1$  for all $ \phi\in \Lambda.$ $(C_b(R),(p_{\phi})_{\phi\in\Lambda} )$ is not m-convex. Let $f\in C_b(R)$ and $\phi\in \Lambda $ such that $ p_{\phi}\left(f\right)=1.$  Since $\phi$ vanishes at infinity, there exists $ x_0\in R$ such that $\vert f\left(x_0\right)\phi\left(x_0\right)\vert=1. $ Show that $\vert f(x_{0})\vert=\sup \lbrace\vert g(x_{0})\vert, p_{\phi}(g)\leq 1\rbrace.$ The proof is due to Beddaa [5], it is given here for completeness. We have $\vert f\left(x_0\right)\vert\leq\sup \lbrace\vert g(x_{0})\vert, p_{\phi}(g)\leq 1\rbrace$  since $p_{\phi}\left(f\right)=1.$  If $p_{\phi}\left(g\right) \le 1=\vert f\left(x_0\right)\phi\left(x_0\right)\vert$, then $ \vert g\left(x_0\right)\phi\left(x_0\right)\vert\le \vert f\left(x_0\right)\phi\left(x_0\right)\vert,$ so $\vert g\left(x_0\right)\vert\le \vert f\left(x_0\right)\vert$ since $\phi\left(x_0\right)>0.$ We have shown that $\vert\delta_{x_0}\left(f\right)\vert= \sup \lbrace\vert \delta_{x_{0}}(g)\vert, p_{\phi}(g)\leq 1\rbrace$,  where $\delta_{x_0}: C_b\left(R\right)\to C, \delta_{x_0}\left(h\right)=h\left(x_0\right),$  is a nonzero continuous multiplicative linear functional on $C_b\left(R\right).$ Thus $(C_b(R),  (p_{\phi})_{\phi\in\Lambda} )$  is saturated.

\noindent\textbf{•}

Oudadess [2] has introduced the following definition:

\noindent\textbf{•}
 
\noindent\textbf{Definition 2.3 [2].} $E$ is v-saturated if for all $\alpha \in \Lambda $  and $x\in E $ such that  $p_{\alpha }\left(x\right)=1,$ there exist $ m_{0} , m$ in $M^{\ast}=M^{\ast}(E)$ such that $ \vert m_{0}(x)\vert=\sup \lbrace\vert m(y)\vert, p_{\alpha}(y)\leq 1\rbrace. $

\noindent\textbf{•}

We have three remarks about this definition.

\noindent\textbf{}

1.  If $E$ is v-saturated, for $\alpha \in \Lambda $  and $x\in E $ such that  $p_{\alpha}(x)=1,$ there exist $ m_{0}, m $ in $ M^{\ast} $ such that $ \sup \lbrace\vert m(x)\vert, p_{\alpha}(y)\leq 1\rbrace=\vert m_{0}(x)\vert <\infty, $  then $M$ is nonempty since $m\in M$, so the fact of replacing $M$ by $M^{\ast}$ is not justified.

2.  Oudadess [2] did not use anywhere the sup-property of his definition. He only used the following deduced property:  for all $\alpha \in \Lambda $  and $x\in E $ such that  $p_{\alpha }\left(x\right)=1,$ there exists $ m_{0}\in M^{\ast}$ such that $\vert m_{0}(x)\vert\ge 1.$

3.  In [2, Theorem 5.4], Oudadess claimed that if $E$ is complete and v-saturated, then the Gelfand map $G$ from $(E,\left(p_{\alpha }\right)_{\alpha \in \Lambda })$  onto $(\hat{E}, (\hat{p}_{\alpha })_{\alpha \in \Lambda })$ is an algebraic and topological isomorphism. For the proof, he used Theorem 3.5 of [1], but this theorem allows us only to conclude that $G$ is a continuous algebraic isomorphism.

\noindent\textbf{•}

\noindent\textbf{REFERENCES} 
 
\noindent\textbf{•}

\noindent\textbf{•}[1] A. C. Cochran, Representation of A-convex algebras, Proc. Amer. Math. Soc.,vol. 41, no. 2, pp. 473-479, 1973.

\noindent\textbf{•}[2] M. Oudadess, v-saturated uniformly A-convex algebras, Math. Japonica, vol. 35, pp. 615-620, 1991.

\noindent\textbf{•}[3] A. Mallios, Topological algebras, Selected Topics, North-Holland, Amsterdam, 1986.

\noindent\textbf{•}[4] E. A. Michael, Locally multiplicatively-convex topological algebras, Mem. Amer. Math. Soc., no. 11, 1952.
 
\noindent\textbf{•}[5] A. Beddaa, Th\'{e}or\`{e}mes du type Gelfand-Naimark dans des alg\`{e}bres topologiques ou bornologiques, Th\`{e}se de troisi\`{e}me cycle, ENS-Rabat, 1987.

\end{document}